\documentclass[english,12pt]{amsart}
\usepackage{babel,graphicx}
\newtheorem{thm}{Theorem}[section]

\newtheorem{lem}[thm]{Lemma}
\newtheorem{prop}[thm]{Proposition}
\theoremstyle{definition}
\newtheorem{defn}[thm]{Definition}
\theoremstyle{remark}
\newtheorem{rem}[thm]{Remark}

\numberwithin{equation}{section}

\begin{document}

\title[Non-rigidity of hyperbolic surfaces lamination]
{Non-rigidity of hyperbolic surfaces lamination}%
\author{Bertrand Deroin}%
\address{\newline
Unit\'e de Math\'ematiques Pures et Appliqu\'ees \newline
\'Ecole Normale Sup\'erieure de Lyon \newline
46, all\'ee d'Italie \newline
F-69364 Lyon Cedex 07 - France} 
\email{bderoin@umpa.ens-lyon.fr}
\subjclass{57R30; 30F}%
\keywords{Teichm\"uller space, Riemann surface lamination}%

\thanks{The author acknowledges support from the 
{\it Swiss National Science Foundation}. }

\begin{abstract} In this note we prove infinite dimensionality of the Teichm\"uller space of a hyperbolic Riemann surface 
lamination of a compact space having a simply connected leaf.\end{abstract}

\maketitle
\section{Introduction}
A \textit{Riemann surfaces lamination} $\mathcal L$ of a topological space $X$ is an atlas of homeomorphisms 
from open sets covering $X$ to the product of the unit disc by a topological space, called \textit{transverse space}, 
such that the changes of coordinates 
(i) preserve the local fibration by discs and (ii) are holomorphic along the fibers. D.~Sullivan discovers some interesting relations between 
$C^1$-conjugacy classes of dynamical systems of dimension $1$ and the space of conformal structures on a lamination by smooth surface. 
%
%
%
%
In \cite{Su1}, he made the systematic study of the \textit{Teichm\"uller space} of a Riemann surface lamination. It is defined as the space of 
all continuous conformal structures along the leaves (transversally continuous Beltrami coefficients), modulo the group of quasi-conformal 
isotopies tangent to the leaves. It is then proved that if the leaves are \textit{hyperbolic} and the space is \textit{compact}, 
there is a ``laminated" Bers embedding theorem. 
Namely, the natural topology on the Teichm\"uller space is Hausdorff, and it is biholomorphic 
to a non empty open subset of the Banach space of holomorphic quadratic differentials along the leaves and continuous on the total space. 

In this note we prove the infinite dimensionality of Teichm\"uller space for a wide class of Riemann surface laminations of a compact space 
whose leaves are hyperbolic, by producing an infinite dimensional space of continuous holomorphic quadratic differentials. 
\begin{thm} \label{main theorem} Let $\mathcal L$ be a Riemann surfaces lamination of a compact space $X$ with hyperbolic leaves and suppose that 
$\mathcal L$ possesses a simply connected leaf. Then the Banach space of continuous holomorphic quadratic differentials is infinite dimensional. \end{thm}
In fact we prove the theorem under a weaker hypothesis. We prove that there exists a universal constant $C$ such that if 
the hyperbolic Riemann surfaces lamination has a leaf $L$ and a point $x$ on $L$
such that $L$ is not transversally isolated and the radius of injectivity at $x$ for the hyperbolic metric in $L$ is greater than $C$ then 
the theorem is true. 

In his famous memoir \cite{Po1}, H.~Poincar\'e proved the existence of non trivial \textit{meromorphic} quadratic differentials
on a compact Riemann surface of genus greater than $2$. Let us recall what was his original method: a compact Riemann surface of genus at least $2$ is the 
quotient of the unit disc ${\bf D}$ by a discrete cocompact group $\Gamma$ of biholomorphisms. Let $f$ be a meromorphic function on the disc 
which is bounded at infinity and consider the following \textit{fuchsian serie}:
\[ \sigma = \sum_{\gamma \in \Gamma} \gamma ^* (f(z)dz^2).\]
Modulo convergence it defines a meromorphic quadratic differential having poles at the projection of the poles of $f$ on $\Sigma$. H.~Poincar\'e observed 
that the serie $\sigma$ converges: if $z$ is a point of the disc which is not in an orbit of a pole of $f$, the module of 
$f(\gamma(z)) \gamma'(z)^2$ is approximately the \textit{euclidian} area of the image by $\gamma$ of a fixed small disc centered at $z$. So because 
$\Gamma$ is discrete and the disc has finite euclidian area, the serie $\sigma$ converges on compact subsets of the disc. 

In \cite{Gh1}, \'E.~Ghys generalizes this method on a hyperbolic Riemann surface lamination $\mathcal L$ of a compact space having a total transversal $\mathcal T$. 
This enables him to construct continuous meromorphic quadratic differential with given poles on $\mathcal T$. 
He considers on each leaf of the lamination a fuchsian serie \`a la Poincar\'e. Using 
the uniformization theory of A.~Candel~\cite{Ca1} and A.~Verjovsky~\cite{Ve1}, he proves that the meromorphic quadratic differential is continuous 
on the total space, if the family of poles is continuous. 
In this note we essentially follow the same way, the originality being that we construct \textit{holomorphic} quadratic differentials. 
The existence of a total transversal 
is then a consequence of our construction. The difficult part is to show that the fuchsian serie $\sigma$ does not vanish identically. For instance,
in Poincar\'e's proof,
the dimension of the space of quadratic differentials on $\Sigma$ is $3g-3$: there exist many polynomials $f$ such that the serie $\sigma$ vanishes identically.  

The paper is organized as follows. In section \ref{Fuchsian series} we define the fuchsian serie $\sigma(q)$ associated to a bounded 
quadratic differential $q$ on the hyperbolic disc, whose support is $r$-separated. For that we extend each value $q_t$ of $q$ by the holomorphic quadratic 
differential $\tilde{q_t}$ of minimizing $L^1$ norm and we define 
\[ \sigma(q):= \sum_{t\in {\bf D}} \tilde{q_t}.\]
Uniform convergence of $\sigma(q)$ on compact subsets of the disc is proved in lemma \ref{convergence of fuchsian series}. 
In section \ref{continuity of sigma} we present some continuity properties of $\sigma$. We introduce a topology on the space 
$\mathcal{DQ}$ of quadratic differential on the disc bounded by $C>0$ and whose support is $r$-separated, for given constants $r,C>0$.
Then we consider the map $\sigma: \mathcal{DQ} \rightarrow HQ(B)$, where $HQ(B)$ is the Banach space of bounded holomorphic quadratic 
differential on the unit hyperbolic ball $B$, and we prove that it is continuous. In section \ref{laminated 
fuchsian series}, we construct fuchsian series on a hyperbolic Riemann surface lamination of a compact space. Each leaf intersects a transversal  
on a separated subset for the hyperbolic metric.
If $q$ is a continuous quadratic differential on the transversal we define on the universal cover of each leaf the fuchsian serie $\sigma(q)$.
These series are invariant by the fundamental group of the leaves and produce a holomorphic quadratic differential on $X$.
We prove the continuity of $\sigma(q)$ in lemma 
\ref{continuity of laminated fuchsian series} if $q$ vanishes on the boundary of the transversal in the transverse topology. 
We use the unifomization theory of A.~Candel and A.~Verjovsky. This is the crucial part of the paper. 
We end in section \ref{proof of the theorem} by the proof of the theorem. \\

I thank \'Etienne Ghys and Alexei Glutsyuk for many helpful conversations.

\section{Fuchsian series} \label{Fuchsian series}
%
%
Let ${\bf D}=\{z\in {\bf C}/ |z|<1\}$ be the \textit{unit disc} and $\tau$ a holomorphic quadratic differential on ${\bf D}$. The 
\textit{$L^1$-norm} of $\tau$ is the volume of $({\bf D}, |\tau|)$ where $|\tau|$ is considered as a singular metric on ${\bf D}$. While 
${\bf C}$ does not carry any non vanishing $L^1$-holomorphic quadratic differential, the disc does. For instance the norm of the quadratic differential 
$\tau = dz^2$ on the disc is $\pi$. 

If $q$ is a quadratic differential at a point $x$ of the disc, we define $\tilde{q}$ to be the holomorphic quadratic differential extending $q$ 
to ${\bf D}$ and \textit{minimizing} the $L^1$-norm. It exists and is unique because of the strict convexity of the $L^1$-norm 
on the space of holomorphic quadratic differentials. One easily checks that $dz^2$ is the minimizing holomorphic quadratic differential 
extending its value at the point $0$. The other minimizing holomorphic quadratic differentials can be deduced from it 
by using the group of automorphisms of the disc. 

Consider the \textit{hyperbolic metric} 
\[ g= \frac{|dz|^2}{(1-|z|^2)^2}\]
on the disc, and $T$ a $r$-separated subset of $({\bf D},g)$, meaning that two distinct points of $T$ are separated by a distance greater than $r$. 
\begin{defn}\label{fuchsian series-def} Let $q=\{q_t\}_{t\in T}$ be a family of quadratic differentials on $T$, 
we define the \textit{Fuchsian serie} $\sigma(q)$ by 
\[ \sigma(q) := \sum_{t\in T} \widetilde{q_t}.\]\end{defn}  
The following lemma is the main result of this section, and can be attributed to H.~Poincar\'e~\cite{Po1}. 
We note $|q|_{\infty,A}$ the uniform hyperbolic norm of a quadratic differential $q$
defined on a subset $A\subset {\bf D}$.
\begin{lem}[Poincar\'e] If $q$ is bounded in hyperbolic norm, 
then $\sigma(q)$ converges uniformly on compact subsets of the disc to a holomorphic quadratic differential 
verifying 
\[ |\sigma(q)|_{\infty,{\bf D}} \leq C(r) |q|_{\infty,T}\]
and 
\[ \label{sigma approximating q} |\sigma(q)-q|_{\infty,T} \leq D(r) |q|_{\infty,T}\]
where $C(r)$ and $D(r)$ are constants depending only on $r$ and $D(r)$ tends to $0$ when $r$ goes to infinity.\label{convergence of fuchsian series}\end{lem}
\noindent\textit{Proof.} First we check some basic properties of the minimizing holomorphic quadratic differentials:\\

(i)  Let $q$ be a quadratic differential at a point $x$ of the disc. The hyperbolic norm of
$\tilde{q}$ is a function depending only on the distance to $x$. It can be seen by using the group of isometries of the disc fixing $x$. 
In particular, if we consider two (non vanishing) quadratic differentials $q_1$ and $q_2$ at points $x_1$ and $x_2$ respectively, 
then we have the symmetry property 
\[ \frac{|\widetilde{q_1}(x_2)|}{|q_1|} = \frac{|\widetilde{q_2}(x_1)|}{|q_2|}.\]\\

(ii)  The hyperbolic norm of a $L^1$ quadratic differential $\tau$ on the disc is an integrable function with respect to the hyperbolic volume $v_g$
because of the formula 
\[ \int_{{\bf D}} |\tau| dv_g = |\tau |_1.\]
In particular, using the minimizing quadratic differential $\tau=dz^2$ of $dz^2$ at the point $0$, we find the general formula 
\[ \int_{{\bf D}} |\tilde q| dv_g = \pi |q|\]
for any quadratic differential $q$ defined at some point $x$. \\

(iii)  By Cauchy formula, there exists a constant $A(r)$ depending only on $r$ such that 
\[ |\tau (0)| \leq A(r) \int _{B_g(0,r)} |\tau |dv_g,\]
for every holomorphic quadratic differential $\tau$ defined on the hyperbolic ball $B_g(0,r)$ of radius $r$ centered at $0$. Again by homogeneity 
of the hyperbolic disc we have the uniform inequality 
\[ |\tau(x)| \leq A(r) \int _{B_g(x,r)} |\tau |dv_g,\]
for every holomorphic quadratic differential $\tau$ defined on the hyperbolic ball $B_g(x,r)$ of radius $r$ centered at $x$. We choose the function 
$A$ to be a decreasing function of $r$.\\

The lemma follows easily from these properties: let $x\in {\bf D}$ and $q$ a quadratic differential of norm $1$ at $x$. 
By (i), 
\[ \sum_{t\in T} |\widetilde{q_t}(x)| \leq |q|_{\infty,T} \sum _{t\in T} |\tilde{q}(t)|,\]
and by (iii), we obtain 
\[\sum_{t\in T} |\widetilde{q_t}(x)| \leq |q|_{\infty,T} A(r/2) \sum _{t\in T} \int_{B_g(t,r/2)} |\tilde q|dv_g.\]
Because $T$ is $r$-separated the balls $B_g(t,r/2)$ are disjoint. Using (ii), we thus have
\[\sum_{t\in T} |\widetilde{q_t}(x)| \leq C(r) |q|_{\infty,T}\]
where $C(r)=\pi A(r/2)$. From this inequality we obtain the uniform convergence of $\sigma(q)$ on compact subsets of the disc. 

Now suppose $x$ is a point of $T$. We have
\[ |\sigma(q)(x) -q(x)|\leq \sum_{t\in T, t\neq x} |\widetilde{q_t}(x)|\]
and doing the same estimates we get 
\[ |\sigma(q)(x)-q(x)|\leq A(r/2) |q|_{\infty,T} \sum_{t\in T,t\neq x} \int_{B_g(t,r/2)} |\tilde q|dv_g.\]
Because $T-\{x\}$ is $r$-separated from $x$, we obtain 
\[ |\sigma(q)(x)-q(x)|\leq A(r/2) |q|_{\infty,T} \int_{{\bf D}-B_g(x,r/2) } |\tilde q|dv_g.\]
But by homogeneity of the disc, $B(R)=\int_{{\bf D}-B_g(x,R) } |\tilde q|dv_g$ does not depend on $x$ and because $\tilde q$ is integrable 
it tends to $0$ when $R$ goes to infinity. We get \ref{sigma approximating q}
with $D(r)=A(r/2)B(r/2)$. But because $A$ is decreasing, $D$ tends to $0$ when $r$ goes to infinity. This ends the proof of the lemma.\\

We conclude the section by proving the following extension result.
\begin{prop}\label{extension on the disc} There exists a constant $r$ such that if $T\subset ({\bf D},g)$ is $r$-separated, 
every bounded quadratic differential 
$q$ on $T$ can be extended to the whole disc to a bounded holomorphic one with an estimate
\[ |q|_{\infty,{\bf D}} \leq E(r) |q|_{\infty,T}.\]\end{prop} 
\noindent \textit{Proof.} Choose $r$ such that $D(r)<1$. Consider the Banach spaces $Q(T)$ of bounded quadratic differential on $T$ 
and $HQ({\bf D})$ of bounded holomorphic quadratic differentials on the discs, both equipped with the $L^{\infty}$ norm. If 
$i: HQ({\bf D})\rightarrow Q(T)$ is the restriction map and $\sigma : Q(T)\rightarrow HQ({\bf D})$ is the map constructed in 
\ref{convergence of fuchsian series} then we have 
\[ ||i\circ \sigma -id|| \leq D(r),\]
where $||.||$ is the norm of operator. So because $D(r)<1$ the operator $i\circ \sigma$ is invertible and $i$ is thus surjective. The proposition is proved.    
\section{A technical lemma}\label{continuity of sigma}
Let $r>0$ and $C>0$ be given constants and define $\mathcal{DQ}$ to be the space of quadratic differentials on the disc, whose hyperbolic norm is bounded 
everywhere by $C$ and whose support is $r$-separated. In this section we construct a \textit{topology} on $\mathcal{DQ}$, and we prove that the map 
$\sigma$ is continuous with respect to it.

Let $K$ be a compact subset of the disc,
$0<\alpha< r/2$ a real number and $q$ an element of $\mathcal{DQ}$. Define the subset $U_{K,\alpha,q}\subset \mathcal{DQ}$ of elements $p$ carrying 
the following properties:\\

(iv)  For each $t\in \mathrm{supp}(q)\cap K$, there exists at most one element of $\mathrm{supp}(p)$ $\alpha$-close to $t$. In the case 
$B_g(t,\alpha)\cap \mathrm{supp}(p)$ is empty, then $|q_t|<\alpha$ 
and if $B_g(t,\alpha)\cap \mathrm{supp}(p)=\{s\}$, then $d(q_t,p_s)<\alpha$. \\

(v)   For each $s\in \mathrm{supp}(p)\cap K$, there exists at most one element of $\mathrm{supp}(q)$ $\alpha$-close to $s$. 
In the case $B_g(s,\alpha)\cap \mathrm{supp}(q)$ is empty, then $|p_s|<\alpha$ 
and if $B_g(s,\alpha)\cap \mathrm{supp}(q)=\{t\}$, then $d(q_t,p_{s})<\alpha$. \\

The family of sets $U_{K,\alpha,q}$ defines a base of open sets for our topology.
The distance $d$ on the total space of the fiber bundle $T^*{\bf D}^{\otimes 2}$ that we consider has not a real significance for our purpose, 
because $K$ is compact. Thus we will work with the following one:
\[ d(q_1,q_2):= |\tilde{q_1}-\tilde{q_2}|_{L^1},\]
for every quadratic differentials $q_1$ and $q_2$ given at some points $x_1$ and $x_2$ respectively. 

Let $B=B_g(0,1)$ be the hyperbolic unit ball centered at $0$, and $HQ(B)$ 
be the space of holomorphic quadratic differentials on $B$ with bounded hyperbolic norm. We equip $HQ(B)$ with the uniform topology. 
\begin{lem}\label{continuity of fuchsian series} The map $\sigma : \mathcal{DQ} \rightarrow HQ(B)$ is continuous.\end{lem}
\noindent \textit{Proof.} Let $q\in \mathcal{DQ}$, $R>r$ and $0<\alpha<r/2$ be real numbers. 
Now choose $p$ in the set $U_{\overline{B}_g(0,R),\alpha,q}$. We will show that $|\sigma(p)-\sigma(q)|_{\infty,B}$ tends to $0$ uniformly 
when $R$ tends to infinity and when $\alpha = o( 1/ \mathrm{vol}(B_g(0,R)) )$. This suffices to prove the lemma. 

Define $T$ to be in $\overline{B}_g(0,R)$ the union of the points of $\mathrm{supp}(q)$ and of the points of $\mathrm{supp}(p)$ 
from which there is no $\alpha$-close point of $\mathrm{supp}(q)$. Similarly define $S$ to be in $\overline{B}_g(0,R+\alpha)$ the union of the 
points of $\mathrm{supp}(p)$ and of the points of $\mathrm{supp}(q)$ from which there is no $\alpha$-close point of $\mathrm{supp}(p)$. 
By (iv) we can thus associate to every point $t$ of $T$ a unique point $s(t)$ of $S\cap B_g(t,\alpha)$. The map $s:T\rightarrow S$
is injective and by (v) its image contains $S\cap \overline{B}_g(0,R-\alpha)$. We thus have, for each $x$ in the unit hyperbolic ball 
$B$ 
\begin{equation}\label{eq 1} 
|\sigma(q)(x)-\sigma(p)(x)|\leq 
\end{equation} 
\[\sum _{t\in T} |\widetilde{q_t}(x) -\widetilde{p_{s(t)}}(x)| 
+\sum_{t\in \mathrm{supp}(q)-B_g(0,R)} |\widetilde{q_t}(x)|+\sum_{s\in \mathrm{supp}(p)-B_g(0,R-\alpha)} |\widetilde{p_s}(x)|.\]

We first bound the first term of the right hand side of \ref{eq 1}. By definition of the map $t\in T\mapsto s(t)\in S$, and by property (iv), 
\[ |\widetilde{q_t}- \widetilde{p_{s(t)}}|_1 <\alpha,\]
so that using the uniform Garding inequality (iii), we find 
\[ |\widetilde{q_t}(x)- \widetilde{p_{s(t)}}(x)|\leq A(\infty) \alpha.\]
In order to bound on the number $|T|$ of points of $T$, observe that $T$ is a priori only a $\alpha$-separated subset, but by definition it is the union of 
two $r$-separated subsets of the ball $B_g(0,R)$. Thus we obtain
\[ \sum _{t\in T} |\widetilde{q_t}(x) -\widetilde{p_{s(t)}}(x)|\leq A(\infty)\alpha |T|\leq 2A(\infty)\alpha 
\frac{\mathrm{vol}(B_g(0,R+r/2))}{\mathrm{vol}(B_g(0,r/2))}.\]

For the other terms of \ref{eq 1}, we proceed like in lemma \ref{convergence of fuchsian series} and we get 
\[ \sum_{t\in \mathrm{supp}(q)-B_g(0,R)} |\widetilde{q_t}(x)|\leq A(r/2)B_g(R-\alpha -1) |q|_{\infty} ,\]
and 
\[ \sum_{s\in \mathrm{supp}(p)-B_g(0,R-\alpha)} |\widetilde{p_s}(x)|\leq A(r/2)B_g(R-2\alpha -1) |q|_{\infty}. \]
Recall that the elements of $\mathcal {DQ} $ are supposed to be bounded by a universal constant $C>0$ so that we conclude
\[ |\sigma(q)-\sigma(p)|_{\infty,B}\leq F(o(R) +\alpha \mathrm{vol}(B_g(0,R))),\]
for a constant $F$ depending only on $r$ and $C$. The lemma is proved.

\section{Continuity of laminated fuchsian series}\label{laminated fuchsian series}
Let $\mathcal L$ be a Riemann surfaces lamination of a topological space $X$.
\begin{defn} A \textit{transversal} of $\mathcal L$ is a closed subset $\mathcal T \subset X$ satisfying the following property: 
there exists a cover of $X$ by flow boxes ${\bf D}\times T$ such that $\mathcal T$ intersects each slice ${\bf D}\times \{t\}$ 
in at most one point. \end{defn} 
On a transversal there is the natural \textit{transverse topology}. For instance, a point $t$ of $\mathcal T$ is a 
boundary point (resp. interior point) if there exists a flow box ${\bf D}\times T$ containing $t$ and such that 
$t$ projects on a boundary point (resp. interior point) of the image of $\mathcal T$ by the projection ${\bf D}\times T\rightarrow T$. 
The boundary of $\mathcal T$ in the transverse topology is denoted by $\partial_t \mathcal T$, and the interior $\mathrm{Int}_t (T)$. 

Suppose now that the leaves of $\mathcal L$ are \textit{hyperbolic} and the space $X$ is \textit{compact}. A.~Candel~\cite{Ca1} and A.~Verjovsky~\cite{Ve1} 
have shown that 
the hyperbolic metric on the leaves is continuous. Thus every leaf with its hyperbolic metric intersects a transversal $\mathcal T$ in a 
$r$-separated subset, where $r=r(\mathcal T)$ is a positive number depending only on $\mathcal T$. If $q$ is a \textit{bounded} 
quadratic differential on $\mathcal T$, we consider on the universal cover $\pi: \widetilde{L}\rightarrow L$ of each leaf the serie 
\[ \sigma(q) = \sum_{t\in \pi ^{-1}(L\cap \mathcal T)} \widetilde{q_{\pi(t)}}.\]
Because this serie converges (Lemma \ref{convergence of fuchsian series}) and is invariant by the fundamental group of $L$,
it induces a holomorphic quadratic differential the leaves of $\mathcal L$. Here is the main result of the section. 
\begin{lem} \label{continuity of laminated fuchsian series} Suppose $q$ is a continuous quadratic differential defined on a 
transversal $\mathcal T$ and vanishing on $\partial _t \mathcal T$. Then $\sigma(q)$ is continuous. \end{lem}
\noindent \textit{Proof.} To prove this lemma, we will use the uniformization theory of A.~Candel and A.~Verjovsky for hyperbolic 
laminations of a compact space. Take a point $x\in X$ and a transversal $T$ containing $x$ in its interior for the transverse 
topology. If we choose a continuous family $v_t$ of unitary vectors of the unitary tangent bundle $U\mathcal L$ of $\mathcal L$ on $T$, 
then we can consider the map 
\[ \pi :{\bf D}\times T\rightarrow X\]
defined on each plaque ${\bf D}\times \{t\}$ as being the covering sending $(0,t)$ to $t$ and the vector 
$\frac{\partial}{\partial x}\in T_0{\bf D}$ of the unit disc to $v_t$. The uniformization result of A.~Candel and A.~Verjovsky \cite{Ca1,Ve1} 
is the fact that this map 
induces by restriction to ${\bf D}\times \mathrm{Int}_t (T)$ a local biholomorphism of Riemann surfaces laminations. In particular,
if $\varepsilon >0$ is chosen small enough, the restriction of $\pi$ to $B_g(0,\varepsilon) \times \mathrm{Int}_t (T)$, where $B_g(0,\varepsilon)$ 
is the hyperbolic ball of radius $\varepsilon$, defines a coordinate 
in which we will prove that $\sigma(q)$ is continuous and the lemma will follow. 

First the preimage $\hat{\mathcal T}$ of $\mathcal T$ by $\pi$ is a transversal of the trivial lamination ${\bf D}\times T$, and the differential 
$\pi^* q$ is continuous on it and vanishes on $\partial_t \hat{\mathcal T}$. Furthermore, $q$ is bounded (let say by $C$), and the intersection 
of $\mathcal T$ with each leaf is $r$-separated, so that the family of quadratic differentials $\pi^*q_{|{\bf D}\times \{t\}}$ defines a map from $T$ to 
the space $\mathcal DQ$ that we have considered 
in the section \ref{continuity of sigma}. Because of Lemma \ref{continuity of fuchsian series} it suffices to prove that it is continuous. 

Take $t_0\in \mathrm{Int}_t(T)$, and $R>r$, $0<\alpha<r/2$ two real numbers. The ball $\overline{B}_g(0,R)\times \{t_0\}$ intersects $\hat{\mathcal T}$ 
in a finite number of points, let say $(z_i,t_0)$, $i=1,\ldots,N$. Let $T_i$ be a neighborhood of $t_0$ in $T$ such that the flow box 
$U_i=B_g(z_i,\alpha)\times T_i$ contains at most one point of $\hat{\mathcal T}$ in every plaque $B_g(z_i,\alpha)\times \{t\}$. Because $\pi^*q$ is continuous 
on $\hat{\mathcal T}$, we can also choose $T_i$ in such a way that if $(z,t)$ and $(z',t')$ are two points of $\hat{\mathcal T}\cap U_i$
then 
\[ d(\pi^* q_{(z,t)}, \pi^* q_{(z',t')})<\alpha.\]
Consider now $\mathcal V= T_1\cap \ldots \cap T_N$. The construction of the $T_i$'s and the vanishing of $\pi^*q$ on $\partial _t \hat{\mathcal T}$ 
is exactly saying that for each $t\in \mathcal V$, $\pi^*q_{|{\bf D}\times \{t\}} \in \mathcal U_{\overline{B}_g(0,R),\alpha,q_0}$, $q_0$ being the 
quadratic differential $\pi^*q_{|{\bf D}\times \{t_0\}}$. Thus we have proved continuity 
of the map $t\in T\mapsto \pi^*q_{|{\bf D}\times \{t\}}\in \mathcal{DQ}$ and therefore the continuity of $\sigma(q)$ by 
using Lemma \ref{continuity of fuchsian series}.

\section{Proof of the theorem}\label{proof of the theorem}
The proof of the theorem is a variation of the open image theorem that we used in proposition \ref{extension on the disc}. If $L$ is a simply connected leaf 
of $\mathcal L$, there exists a transversal $T_1$ cutting $L$ in an interior point $t_1$ of $T_1$ with $r=r(T_1)$ arbitrarily large, let say such that 
$D(r)<1$. We can suppose furthermore that the bundle $T^*\mathcal L ^{\otimes 2}$ is trivial on $T_1$ if for instance $T_1$ is chosen in a flow box. 

First suppose that $L$ is not transversally isolated (i.e. $t_1$ is not isolated in $T_1$), and take a closed neighborhood $T_2$ of $t_1$ in $T_1$. 
We will show that every continuous quadratic differential on $T_2$ can be extended to $X$ to a continuous holomorphic one, and this will prove the result 
because $T_2$ is infinite. 

We begin by extending $q$ to a continuous quadratic differential $\overline q$ on $T_1$ vanishing on $\partial _t T_1$ and such that $|\overline q|_{T_1,\infty}
\leq |q|_{T_2,\infty}$. Then we consider 
\begin{equation}\label{first step of approximation of q}\sigma _1= \sigma(\overline q).\end{equation} 
By Lemma \ref{continuity of laminated fuchsian series}, $\sigma _1$ is a continuous holomorphic quadratic differential 
on $X$ such that
\[ |\sigma _1|_{X,\infty}\leq C(r) |q|_{\infty,T_2}\]
and 
\[ |q-(\sigma_1)_{|T_2}|_{\infty,T_2} \leq D(r) |q|_{\infty,T_2}.\]
This is the first step of an approximation process that goes inductively as follows. Given continuous holomorphic quadratic differentials 
$\sigma_1,\ldots, \sigma_n$ on $X$ satisfying for $i=1,\ldots,n-1$
\begin{equation}\label{eq 8} |\sigma _{i+1}|_{X,\infty}\leq C(r) |q - (\sigma_1 +\ldots +\sigma_i)_{|T_2} |_{\infty,T_2}\end{equation}
and 
\begin{equation}\label{eq 9} |q-(\sigma_1+\ldots +\sigma_{i+1})_{|T_2}|_{\infty,T_2} \leq D(r) |q- (\sigma_1 +\ldots +\sigma_i)_{|T_2} |_{\infty,T_2},\end{equation}
we construct $\sigma_{n+1}$ by considering $\sigma_{n+1}=\sigma_1(q- (\sigma_1 +\ldots +\sigma_n)_{T_2})$ rather than $q$ (see \ref{first step of approximation of q}). 
By construction both \ref{eq 8} and \ref{eq 9} 
are also verified for $i=n$. Thus we have for each $n\geq 1$,
\[ |q- (\sigma_1 +\ldots +\sigma_n)_{|T_2}|\leq D(r)^n |q|_{\infty,T_2}\]
and
\[ |\sigma_n|_{\infty, X }\leq C(r)D(r)^{n-1} |q|_{\infty,T_2}.\]
The serie $\sum_{n} \sigma_n$ converges uniformly to a continuous holomorphic quadratic differential on $X$ whose restriction to $T_2$ is $q$. 
The proof of the theorem in that case is achieved.  
\begin{rem} In fact we have thus proved the theorem under the weaker hypothesis that the hyperbolic Riemann surface lamination has a leaf $L$ and a point $x$ on $L$
such that $L$ is not transversally isolated and the radius $r$ of injectivity at $x$ for the hyperbolic metric in $L$ verifies $D(r)<1$.\end{rem}
Now, if there exists on $\mathcal L$ a simply connected leaf $L$ which is transversally isolated, then every holomorphic quadratic differential on $L$ 
whose hyperbolic norm tends to $0$ at infinity can be extended by $0$ to a continuous holomorphic quadratic differential on $X$. So the space of continuous holomorphic 
quadratic differentials contains at least the one defined on $L$ by the formula 
\[ \tau = f(z)dz^2,\]
where $z:L\rightarrow {\bf D}$ is a biholomorphism and $f$ is a holomorphic function on the disc bounded at infinity. This completes the proof of the theorem.

\end{document}